\documentclass[11pt, a4paper]{article}

\usepackage{amsfonts,amssymb,amsbsy,amsmath,amsthm,latexsym,bm}
\usepackage{hyperref,url}
\usepackage{graphicx}

\usepackage[T1]{fontenc}

\usepackage[top=3cm, bottom=3cm, left=3cm, right=3cm]{geometry}

\newtheorem{theorem}{Theorem}[section]

\newtheorem{lemma}[theorem]{Lemma}
\newtheorem{proposition}[theorem]{Proposition}
\newtheorem{corollary}[theorem]{Corollary}
\newtheorem{remark}[theorem]{Remark}

\theoremstyle{definition}
\newtheorem{example}[theorem]{Example}
\usepackage{comment}
\usepackage{cleveref}

\newcommand{\tr}{\textup{tr}}

\DeclareMathOperator{\diag}{diag}
\newcommand{\Diag}{\textup{Diag}}

\DeclareMathOperator{\rank}{rank}

\DeclareMathOperator{\vect}{vec}
\DeclareMathOperator{\st}{s.t.}

\newcommand{\R}{\mathbb{R}}
\newcommand{\Z}{\mathbb{Z}}

\title{ Many faces of integer semidefinite programming }

\author{
Mirjam D\"{u}r\thanks{Department of Mathematics, University of Augsburg, Germany. Email: mirjam.duer@math.uni-augsburg.de}
\quad and \quad
Renata Sotirov\thanks{CentER, Department of Econometrics and OR, Tilburg University, The Netherlands. Email: R.Sotirov@tilburguniversity.edu}
}

\date{}

\begin{document}

\maketitle

\begin{abstract}
\noindent
We first establish several properties of integer positive semidefinite matrices, extending known results for binary and ternary positive semidefinite matrices.
We then show how to reformulate an integer quadratic program (IQP) as an equivalent integer semidefinite program (ISDP).
We distinguish three cases, namely IQPs with two, three, or more than three distinct integer values, and provide several ISDPs for each case.  Each reformulation yields either a binary ISDP, a ternary ISDP, or a rescaled version of one of these.
The results show that ternary SDP formulations for problems with three distinct integer values provide  SDP relaxations with fewer constraints and, in the case of equally spaced integers, also smaller matrix variables than binary-based SDP relaxations.
For quadratic problems whose variables can take more than three distinct integer values, one can choose among three formulations based on the size of their matrix variables, the number of constraints, and the strength of the resulting SDP relaxations.
We also show how to reformulate a mixed-integer QP as a mixed-integer~SDP.
\medskip

\end{abstract}

\section{Introduction}

Integer semidefinite programming can be viewed  as a generalization of integer linear programming
in which the vector of variables is replaced by a positive semidefinite (PSD) matrix variable whose
entries are constrained to be integer.
It has recently been shown that  integer semidefinite programs (ISDPs) admit significantly different formulations for several classes of problems compared to traditional modeling approaches, and also enable the development of new algorithmic techniques, see e.g.,~\cite{DeMeijerSotirov2024,de2022chv}.
It is also well known that optimization problems in which variables take values in the set $\{-1,1\}$, such as the max-cut and max-$k$-sat problems, can be formulated as ISDPs~\cite{GoemansWilliamson1995}.
There exist binary SDP formulations for
binary quadratically constrained quadratic programs and binary quadratic matrix programs \cite{DeMeijerSotirov2024}.
Recently, it was shown how to model ternary quadratic problems as ternary  semidefinite programs (TSDPs)~\cite{deMeijer2026BeyondBinarity}.

It is straightforward to embed an  integer linear program (ILP) into the space of diagonal integer PSD matrices, resulting in an ISDP.
However, it remains unclear how to construct a general ISDP for an integer quadratic program~(IQP).
Here, by a general ISDP, we mean an optimization problem involving integer positive semidefinite (PSD) matrices whose entries can take any finite number of distinct integer values.
This difficulty may be related to the limited number of results available on general integer PSD matrices.
Decomposition results of general  integer PSD matrices were presented in
\cite{Mordell1937}, and more recently in \cite{DeLoera2025}, where it was shown that any integer PSD matrix of order $\leq 5$ has a finite decomposition into  rank-one matrices.
It is also known that integer PDS matrices of order greater than~$6$ can be represented as the sum of rank-one matrices and one sporadic matrix.
Binary and ternary PSD matrices of any order admit finite decompositions into rank-one matrices, see \cite{Letchford2012} and~\cite{DeMeijerSotirov2024}, respectively. We do not refer to these two types of matrices as general integer matrices.
In this paper, we exploit existing and new results on integer PSD matrices to reformulate IQPs as binary or ternary SDPs, or their rescaled versions.

There is vast literature on solving  binary quadratic optimization problems. However, there is not much literature on solving IQPs whose variables can take any integer value. Let us mention some of the existing work~\cite{christoph_2011,buchheim2013exact,Buchheim2018SDPbasedBF,Buchheim2012SemidefiniteRF,deMeijer2026BeyondBinarity}.
General nonconvex mixed-integer programming solvers, such as {COUENNE}~\cite{Belotti2009Branching}, {BARON}~\cite{Sahinidis2016BARON}, and Gurobi~\cite{gurobi}, are also used to solve general IQPs.

\paragraph{Main results and outline.}

We first derive several properties of integer PSD matrices by  extending known results for binary and ternary PSD matrices. We further establish integer semidefinite characterizations of positive semidefinite matrices of prescribed rank, along with several other results.

The results on integer PSD matrices are then exploited to reformulate mixed-integer quadratic programs (MIQPs) as  integer semidefinite program (MISDPs). To the best of our knowledge, we are the first to reformulate general (M)IQPs as (M)ISDPs.
The considered quadratic programs have quadratic objective and linear constraints. However, it is straightforward to extend our approach to problems with quadratic constraints.
We distinguish among IQPs with two, three, and more than three distinct integer values.

In Section~\ref{sect:two discrete values}, we reformulate an IQP whose variables can take two distinct integer values as a binary semidefinite program (BSDP) and as a general integer SDP.
In Section~\ref{sect:three discrete values}, we first show how to reformulate an IQP in which each of the $n$ variables can take three distinct integer values as a binary SDP with a matrix variable of order $2n+1$, using one-hot decoding.
Then, for three distinct non-equally spaced integer values, we derive an ISDP with a ternary matrix variable of order $2n+1$, while for equally spaced values, we derive a (scaled) ternary SDP with a matrix variable of order $n+1$.
In Section~\ref{sect:compareSDPThreeGeneral} (resp.~Section \ref{sect:SDP relax of equally spaced}), we compare the SDP relaxations for the non-equally spaced (resp., equally spaced) case and prove their equivalence.
While both SDP relaxations for non-equally spaced integers have matrix variables of the same order, the one-hot-based SDP relaxation has more constraints. The ternary-based SDP relaxation for equally spaced integers has a matrix variable of order $n+1$, compared to $2n+1$ in the one-hot formulation, and fewer constraints.
Therefore, ternary-based SDP relaxations are more favorable for solving IQPs with three distinct integer values.

Section~\ref{sect:generalISDP} considers an IQP  whose variables can take values from a finite set of more than three integers. We use one-hot decoding, binary expansion, and ternary expansion to reformulate the IQP into two different binary semidefinite programs and one ternary semidefinite program, respectively.
The ISDP  derived from one-hot decoding has the largest matrix variable and yields the strongest of the corresponding three SDP relaxations.
However, that relaxation may become too computationally expensive to solve as $n$ or the number of distinct integer values increases. The other two SDP relaxations are incomparable, see Section~~\ref{sect:comparision of general formularins}. In Section~\ref{sect: mixed-integer}, we derive an MISDP formulation of a MIQP.
The complexity of the (M)ISDPs derived here, measured by the size of their integer matrix variables,
grows with the number of distinct integer values and depends on the encoding used, for fixed $n$.

\subsection{Notation}
We use $e_k$ to denote the all ones vector of length $k$.
The identity matrix and the matrix of all ones of order $n$
are denoted by $I_n$ and $J_n,$ respectively.
We omit the subscripts of these matrices when there is no confusion about the order.
We define $E_{ij} :=u_iu_j^T + u_ju_i^T$, where the  i$th$ unit vector is denoted by $u_i$.

The set of all $n \times n$ integer (resp.~real)  symmetric matrices is denoted by  $\mathcal{S}^n(\Z)$ (resp.~$\mathcal{S}^n(\R)$).
The cone of symmetric positive semidefinite matrices is defined as $\mathcal{S}^n_+(\R) := \{ {X} \in \mathcal{S}^n(\R) \, : \, \, {X} \succeq {0} \}$.
The set of integer PSD matrices is denoted by $\mathcal{S}^n_+(\Z)$.
For a matrix $X\in \mathcal{S}^n(\Z)$, we say that $X$ is PSD if and only if it is positive semidefinite when regarded as a real matrix. The trace of a square matrix ${X}=(X_{ij})$ is given by $\tr({X})=\sum_{i}X_{ii}$.
 For any ${X},{Y} \in \mathbb{R}^{n \times m}$ the trace inner product is defined as $\langle {X}, {Y} \rangle := \tr({X}^T {Y})$.  For $X\in \R^{n\times m}$ we denote by $|X|$ the matrix obtained from $P$ by taking the absolute value of each entry.

 We regularly work with  augmented matrices of size $(n+1)\times (n+1)$, in which case the first row and column are indexed by 0.  We define  $[n]:=\{1,\ldots, n\}$.

The operator $\diag : \mathbb{R}^{n \times n} \rightarrow \mathbb{R}^n$ maps a square matrix to a vector consisting of its diagonal elements. We denote by $\Diag : \mathbb{R}^n \rightarrow \mathbb{R}^{n \times n}$ its adjoint operator.
The $\vect (\cdot)$ operator maps an $n\times m$ matrix to a vector of length  $nm$ by stacking the columns of the matrix on top of each other. The Hadamard product of two matrices $X=(x_{ij})$
and $Y=(y_{ij})$ of the same size is denoted by~$\circ$, and defined as follows  $(X\circ Y)_{ij}:= x_{ij}y_{ij}.$
The Kronecker product $A \otimes B$ of matrices $A \in {\R}^{p \times q}$ and $B\in {\R}^{r\times s}$ is defined as the $pr \times qs$ matrix composed of $pq$ blocks of size $r\times s$, with block $ij$ given by $a_{ij}B$, $i \in [p]$, $j \in [q]$.

\section{Preliminaries}

A mixed-integer semidefinite program  is an optimization problem of the following form:
\begin{align} \label{MISDP}
\begin{aligned} \min \quad &  \langle {C}, X \rangle \\
\text{s.t.} \quad &  \langle A_i, X \rangle = b_i \quad i\in [m] \\[0.7ex]
&  {X}_{ij} \in \mathbb{Z} \quad \text{for all } (i,j) \in \mathcal{J} \\[0.7ex]
&  X \succeq {0}, \\
 \end{aligned}
\end{align}
where ${C}, A_i \in {\mathcal{S}^n(\R)}$, ${b} \in \mathbb{R}^m$, $\mathcal{J} \subseteq [n] \times [n]$.
If all variables are restricted to be integer, then the corresponding optimization problem is known as an integer semidefinite program.
If all variables are binary (resp.~ternary) then we refer to the corresponding optimization problem as
binary (resp.~ternary) semidefinite program.
In a ternary SDP, we additionally allow  a nonconvex constraint requiring the diagonal of the matrix variable to coincide with the absolute value of its first row and column \cite{deMeijer2026BeyondBinarity}.

In this section we list known properties  of integer  PSD matrices and derive several new ones.
The integer rank-one decomposition of PSD matrices was already studied in 1937 by Mordell~\cite{Mordell1937}, and further analyzed by De Loera~\cite{DeLoera2025}.
\begin{lemma}[\cite{Mordell1937}]
 If $n\leq 5$, then for any $X \in  \mathcal{S}_+^n(\Z)$, one can find a finite index set $K$ and vectors $x_i \in \mathbb{Z}$, $i\in K$ such that
$X = \sum_{i \in K} x_i x_i^T.$
\end{lemma}
The result of the previous lemma fails for $n=6$ due to the existence of sporadic matrices.  An integer matrix  $X \in  \mathcal{S}_+^n(\Z)$ is sporadic if there do not exist $x\in \Z^n \setminus \{0\}$ such that $X-xx^T \succeq 0.$
For $n=6$ the following matrix
$$ M = J_6+I_6 - E_{12}-E_{23} $$
is the unique  sporadic matrix under unimodular equivalence \cite{Mordell1937}.
We say that two matrices $X$ and $Y$ are unimodularly equivalent if $Y=UX = UXU^T$ for some $U\in {\rm GL}(n,\Z):= \{U \in \Z^{n\times n} \,:\, |\det(U)| =1\}$.
\begin{lemma}[\cite{DeLoera2025}]
If $n = 6$, then  for any  $X \in \mathcal{S}_+^n(\Z)$,
$X = \sum_{i \in K} x_i x_i^T + Y,$
for $x_i \in \Z$ and  $Y$ unimodularly equivalent to
$M$, where $K$ is a finite set.
\end{lemma}
De Loera et al.~\cite{DeLoera2025} prove that the cone of positive semidefinite matrices has a finite generating set of integer points,  under the action of a finitely generated group. They also prove the following result.
\begin{proposition}[\cite{DeLoera2025}]
There is no sporadic point in $\mathcal{S}^n_+(\Z)$  if and only if every integer PSD  matrix in  $\mathcal{S}^n_+(\Z)$
has an integer rank-one decomposition.
\end{proposition}
Since every binary \cite{Letchford2012} and ternary \cite{DeMeijerSotirov2024} PSD matrix has an integer rank-one decomposition, the corresponding  sets of matrices do not have sporadic points.
We list below known results on binary and ternary PSD matrices and derive corresponding results for their scaled versions.
\begin{proposition}[\cite{GoemansWilliamson1995}] \label{rank-onepm1}
Let $X \in \{ -1, 1 \}^{n\times n}$ be a symmetric matrix. Then, $X\succeq 0$ if and only if $X = xx^\top$ for some $x \in \{\pm 1\}^n$.
\end{proposition}
A similar result holds for a rescaled matrix. In particular we have the following result.
\begin{corollary}
\label{Lemma:rescallingpm1}
Let $\alpha \in \mathbb{Z}$, $\alpha > 0$, and let
$X \in \{ -\alpha, \alpha \}^{n\times n}$ be a symmetric matrix. Then, $X\succeq 0$ if and only if $X =\alpha xx^\top$ for some $x \in \{\pm 1\}^n$.
\end{corollary}
\begin{proof}
This follows immediately from Proposition~\ref{rank-onepm1}.
\end{proof}
The following result characterizes binary SDP matrices of rank-one.
\begin{proposition}[\cite{DeMeijerSotirov2024}]
\label{prop:Rank1}
Let $Y=\begin{bmatrix}
1 & \diag({X})^T  \\ \diag({X}) & {X}
\end{bmatrix} \succeq {0}$.
Then ${\rank(Y)} = 1$ if and only if ${X} \in \{0,1\}^{n \times n}$.
\end{proposition}
We further generalize the results of Proposition~\ref{prop:Rank1} to rescaled matrices.
\begin{lemma} \label{lemma:alphaS0-1}
Let $\alpha \in \mathbb{Z}$, $\alpha \neq 0$,  and
$X\in {\mathcal S}^n(\mathbb{Z})$ such that  $X_{ij} \in \left \{0, \alpha^2 \right \}^{n\times n}$ for all $i,j\in [n]$. Let
$Y=\begin{bmatrix} 1& x^T  \\  x & {X}
\end{bmatrix}$ where $\diag(X)=\alpha x$.
Then, ${\rank(Y)} = 1$ if and only if $Y\succeq 0$.
\end{lemma}
\begin{proof}
It follows from the assumptions of the lemma that $x\in \{0,\alpha\}.$
Assume $Y\succeq 0$ and $Z\in \{0,1\}^{n\times n}$ such that $X=\alpha^2 Z$.
The Schur complement lemma implies $X - xx^T \succeq 0.$
Then,
$$
\begin{aligned}
X -  xx^T &= \alpha^2 Z - \frac{1}{\alpha^2} \diag(X)\diag(X)^T
= \frac{1}{\alpha^2} ( \alpha^4 Z - \diag(X)\diag(X)^T  ) \\
& = \frac{1}{\alpha^2} ( \alpha^4 Z - \alpha^4 \diag(Z)\diag(Z)^T  )
= \alpha^2 (Z-\diag(Z)\diag(Z)^T ).
\end{aligned}
$$
Since $Z$ is a binary matrix we have  $Z_{ii}=Z_{ii}^2$ for all $i\in [n]$,
from where it follows that
$Z -\diag(Z)\diag(Z)^T =0$ and consequently $\rank(Z)=1$.
Denote $z:=\diag(Z)$.
Then, $X=\alpha^2 Z = \alpha^2 z z^T = xx^T $ where $x=\alpha z$,
from where it follows that $\rank(Y)=1$.

The converse is trivial.
\end{proof}
We derive rank results beyond the rank-one case for integer PSD matrices with entries in $\{0,\alpha^2\}$.
\begin{lemma}
Let  $\alpha \in \mathbb{Z}$, $\alpha\neq 0$, and let $X\in{\mathcal S}^n(\mathbb{Z})$ be such that
$X\in\{0,\alpha^2\}^{n\times n}$.
Let
$ Y= \begin{bmatrix}
r & x^T\\
x & X
\end{bmatrix}$
where  $\diag(X)=\alpha x$, and $r\in \mathbb{N}$, $r\geq 2$.
If $Y\succeq 0$, then
$
X\succeq 0$ and $\rank(X)\le r$.
\end{lemma}
\begin{proof}
The assertion that $X\succeq 0$ is immediate.
Define $Z:=X/\alpha^2$. Then $Z\in\{0,1\}^{n\times n}$ and $Z\succeq 0$.
It follows from $\diag(Z)=\diag(X)/\alpha^2=x/\alpha$ that $x=\alpha\diag(Z)$.

Consider the congruence scaling
$$
\begin{bmatrix}
1 & 0\\
0 & \frac{1}{\alpha}I
\end{bmatrix}
Y
\begin{bmatrix}
1 & 0\\
0 & \frac{1}{\alpha}I
\end{bmatrix}
=
\begin{bmatrix}
r & \diag(Z)^T\\
\diag(Z)& Z
\end{bmatrix}\succeq 0.
$$
It follows from the above linear matrix inequality and the binary rank-$r$ result \cite{DeMeijerSotirov2024}[Corollary 2] that $\rank(Z)\le r.$
Since $X=\alpha^2 Z$, it follows that  $\rank(X)=\rank(Z)$. Therefore
$\rank(X)\le r.$
\end{proof}
The following result extends a similar result from \cite{DeMeijerSotirov2024}[Corollary 3].
\begin{proposition}\label{prop:ExactRank_rescaledBin}
Let  $\alpha \in \mathbb{Z}$, $\alpha\neq 0$ and let $X\in{\mathcal S}^n(\mathbb{Z})$ be such that
$X_{ij}\in\{0,\alpha^2\}^{n \times n}$.
If there exists  a matrix $P\in\{0,\alpha\}^{n\times r}$
such that
$$
\frac{1}{\alpha} P^T e \ge  e,
\quad
\alpha P e = \diag(X)
\quad \mbox{ and } \quad
 Y=
\begin{bmatrix}
I_r & P^T\\
P & X
\end{bmatrix}
\succeq 0,
$$
then $ X\succeq 0$ and $\rank(X)=r.$
\end{proposition}

\begin{proof}
By the Schur complement lemma we have  $X-PP^T \succeq 0,$
from where it follows
$X_{ii}-(PP^T)_{ii} \geq 0$ for every $i\in[n]$.
Since $X_{ii}\in\{0,\alpha^2\}$ and
$(PP^T)_{ii} =\sum_{j=1}^r P_{ij}^2,$ it follows that   every row of $P$ contains at most one nonzero entry. Consequently, the columns of $P$ have pairwise disjoint supports.
Moreover,
$\frac{1}{\alpha} P^T e \ge e$
implies that every column of $P$ is nonzero. Hence the columns of $P$ are nonzero and mutually orthogonal, and therefore
$\rank(P)=r.$

From
$\alpha P e = \diag(X),$ it follows that
$\alpha(P e)_i=X_{ii}$
for all $ i\in[n].$
 Because each row of $P$ contains at most one nonzero entry,
 we have
$ (PP^T)_{ii} = \alpha(P e )_i = X_{ii},$
and therefore $\diag(X-PP^T)=0.$
Since $X-PP^T \succeq0$,  it follows that  $X-PP^T=0,$ and thus $X=PP^T.$
Consequently, $\rank(X) = \rank(PP^T) = \rank(P)=r.$
\end{proof}

The rest of this section provide characterisation of ternary PDS matrices.
\begin{proposition}[\cite{DeMeijerSotirov2024}]
\label{prop:Rank1pm1}
Let $Y = \begin{bmatrix} 1 &  x^T  \\ {x} & {X}\end{bmatrix}$ be a symmetric matrix  with $\diag(X)=|x|$.
Then,
$Y \in \{0,\pm 1\}^{(n+1)\times (n+1)}$, $Y \succeq {0}$ if and only if $X=xx^T.$
\end{proposition}
Next, we derive a rank-one result for the case that PSD integer matrices have elements in $\{-\alpha^2, 0, \alpha^2 \}$ where $\alpha \in \mathbb{Z}.$
\begin{lemma} \label{lemma:alphaSminalpha0-1}
Let $\alpha \in \mathbb{Z}$, $\alpha \neq 0$  and  $X\in {\mathcal S}^n(\mathbb{Z})$ such that
$X_{ij} \in \left \{-\alpha^2, 0, \alpha^2 \right \}^{n\times n}$.
Let $Y=\begin{bmatrix} 1 & x^T  \\  x & {X}\end{bmatrix}$ where $\diag(X)=  |\alpha x|$.
Then, ${\rank(Y)} = 1$ if and only if $Y\succeq 0$.
\end{lemma}
\begin{proof}
Note that the assumptions of the lemma imply that
$x\in \{-\alpha, 0,\alpha\}.$
Assume that  $Y\succeq 0$ and $Z\in \{0,\pm 1\}^{n\times n}$ such that $X=\alpha^2 Z$.
Then by using a similar derivation and reasoning as in Lemma \ref{lemma:alphaS0-1} we
obtain that $\rank(Z)=1$ from where it follows the result.
The converse direction is trivial.
\end{proof}
We now present a rank result for $\{-\alpha^2, 0, \alpha^2 \}$-PSD matrices  analogous to Proposition~\ref{prop:ExactRank_rescaledBin} for $\{ 0, \alpha^2 \}$-PSD matrices.
\begin{proposition}\label{Prop:ExactRank_signed_rescaledTer}
Let $\alpha \in \mathbb{Z}$, $\alpha \neq 0$  and  $X\in {\mathcal S}^n(\mathbb{Z})$ satisfy
$X \in \left \{-\alpha^2, 0, \alpha^2 \right \}^{n\times n}$.
If there exists a matrix $P\in\{0,\alpha\}^{n\times r}$   such that
\[
|P|^T e \ge |\alpha| e,
\qquad
|P|e=\frac{1}{|\alpha|}\diag(X),
\quad \mbox{ and } \quad
Y=
\begin{bmatrix}
I_r & P^\top\\
P & X
\end{bmatrix}
\succeq 0.
\]
Then
$ X\succeq 0$ and $\rank(X)=r.$
\end{proposition}
\begin{proof}
Similarly to the proof of Proposition  \ref{prop:ExactRank_rescaledBin}.
\end{proof}

\section{IQPs with two integer values as ISDPs }
\label{sect:two discrete values}
We reformulate an integer quadratic program with linear constraints having two distinct integer values as a binary semidefinite program.

Let  $\alpha, \beta \in \mathbb{Z}$ such that $\alpha < \beta$, $\alpha \neq 0$, $\beta \neq 0$.  Consider the following IQP:
\begin{equation} \label{QP_with_2_values}
   \begin{array}{cl}
      \min & x^T Q x \\
      \st    & Ax = b \\
              & x \in \{ \alpha, \beta \}^n,
   \end{array}
\end{equation}
where $Q \in {\mathcal S}^n(\mathbb{R})$, $A \in \R^{ m \times n }$, $b \in \R^m$.

Observe that we can map $\{ \alpha, \beta \}^n$ to $\{ 0,1 \}^n$ via the affine mapping $\varphi : \R^n \to \R^n$ with $\varphi(u) = \frac{1}{\beta - \alpha} u - \frac{\alpha}{\beta - \alpha} e$ which has the inverse $\varphi^{-1}(v) = (\beta - \alpha) v + \alpha e$.
Using the variable transformation $x := \varphi^{-1} (y)$, the objective becomes
\[
x^T Q x = (\beta - \alpha)^2 y^T Q y + 2 \alpha (\beta - \alpha) y^T Q e + \alpha^2 e^T Q e \]
and the linear constraints can be written as
$(\beta - \alpha) Ay = b - \alpha Ae.$
Defining $\bar y^T := (1, y)^T \in \{0,1\}^{n+1}$,
we obtain the following binary QP, which is equivalent to~\eqref{QP_with_2_values}:
\begin{equation}  \label{binary_equivalent_of_QP_with_2_values}
   \begin{array}{cl}
      \min & \bar{y}^T \bar{Q} \bar{y} \\[1ex]
      \st    & A y = \frac{1}{\beta- \alpha}(b-\alpha Ae) \\[1ex]
              &  \bar{y}^T = [1,y^T] \in \{ 0,1 \}^{n+1}
   \end{array}
\end{equation}
were
$$ \bar{Q} = \begin{bmatrix}
			\alpha^2 e^T Q e &  \alpha (\beta - \alpha) (Qe)^T \\
			\alpha (\beta - \alpha) Qe & (\beta - \alpha)^2 Q
    \end{bmatrix} \in  {\mathcal S}^{n+1}(\mathbb{R}).
$$

\medskip
The quadratic objective in \eqref{binary_equivalent_of_QP_with_2_values} can be written as
$\langle \bar{Q} , \widetilde{Y} \rangle $ where we  substitute
$\widetilde{Y}$ for $\bar{y}\bar{y}^T$.
Now, by exploiting Proposition~\ref{prop:Rank1} we can rewrite \eqref{binary_equivalent_of_QP_with_2_values} equivalently as the following BSDP:
\begin{equation} \label{BSDP_with_2_values}
   \begin{array}{cl}
      \min & \langle \bar{Q}, \widetilde{Y}\rangle \\[1ex]
      \st    & A y = \frac{1}{\beta- \alpha}(b-\alpha Ae) \\[1ex]
              & \widetilde{Y} = \begin{bmatrix} 1 & y^T \\ y & Y \end{bmatrix} \succeq 0, ~~ \diag (Y) = y \\[2ex]
            & y \in \{ 0,1 \}^{n}.
   \end{array}
\end{equation}
Here, the integrality of the off-diagonal elements in $Y$ follows from the binarity of $y$ and the positive semidefiniteness constraint, see e.g.,~\cite{DeMeijerSotirov2024}[Section 3.1].
Thus, off diagonal elements of  $Y$ are so-called implied integers.
Interested reader can view~\cite{DeMeijerSotirov2024} for binary SDP reformulations of various binary optimization problems.

However, one can also exploit Lemma~\ref{lemma:alphaS0-1} to reformulate IQP \eqref{QP_with_2_values} as an ISDP.
To this end, we introduce the affine  mapping $\rho : \R^n \to \R^n$ that maps $\{ \alpha, \beta\}^n$ to $\{ 0, \beta - \alpha\}^n$, defined by $\rho(u)= u-\alpha e$, whose inverse is $\rho^{-1}(v) = v+\alpha e.$
Proceeding similarly to the previous case,  by rewriting the objective function and the constraints, we obtain the following integer SDP reformulation of the IQP \eqref{QP_with_2_values}:
\begin{equation} \label{ISDP_with_2_values}
   \begin{array}{cl}
      \min & \langle \widetilde{Q}, \widetilde{Y}\rangle \\[1ex]
      \st    & A y = b-\alpha Ae \\[1ex]
              &\widetilde{Y}  = \begin{bmatrix} 1 & y^T \\ y & Y \end{bmatrix} \succeq 0, ~~ \diag (Y) = (\beta -\alpha) y \\[2ex]
            & Y \in \left \{ 0, (\beta - \alpha)^2 \right \}^{n\times n},
   \end{array}
\end{equation}
where
$$
\widetilde{Q} = \begin{bmatrix}
		\alpha^2 e^T Q e &  \alpha e^T Q \\
			\alpha  Qe &  Q
   \end{bmatrix}  \in  {\mathcal S}^{n+1}(\mathbb{R}).
$$
This formulation is the rescaled formulation \eqref{BSDP_with_2_values}.
\begin{remark}
Observe that if either $\alpha=0$ or  $\beta=0$, then, after appropriate scaling, the ISDP reformulation of \eqref{QP_with_2_values} reduces to BSDP~\eqref{BSDP_with_2_values}.
Moreover,  in the case of $\alpha = - \beta$, one can use  Lemma~\ref{Lemma:rescallingpm1} to reformulate the IQP~\eqref{QP_with_2_values} as an ISDP with matrix variables of order~$n$.
\end{remark}

\section{IQPs with three integer values  as ISDPs}
\label{sect:three discrete values}

We reformulate an IQP in which each variable takes one of three distinct integer values, as both a binary and a ternary SDP.
In  Sections \ref{sect:compareSDPThreeGeneral}  and \ref{sect:SDP relax of equally spaced},  we compare the resulting SDP relaxations and show that they are equivalent.

Consider the following IQP in which the variables are constrained to take  one of three distinct integer values, $\alpha$, $\beta$ or $\gamma$.
In particular, consider the following IQP:
\begin{equation} \label{QP_with_3_values}
   \begin{array}{cl}
      \min & x^T Q x \\
      \st    & A x = b \\
              & x \in \{\alpha, \gamma, \beta \}^n,
   \end{array}
\end{equation}
where $Q \in {\mathcal S}^n (\mathbb{R})$, $A \in \R^{ m \times n }$, $b \in \R^m$.
In what follows, we consider several integer SDP reformulations of \eqref{QP_with_3_values}.
While our binary one-hot reformulation applies to any three integer values, our ternary formulations depend on whether these values are equally spaced.

\medskip
\noindent {\bf General binary reformulation.}
We may assume, w.l.o.g.,~that $\gamma=0$.
Let us now introduce two binary vectors $v^1, v^2 \in \{0,1\}^n $ that satisfy $v^1+v^2 \leq e$, and define
\begin{align} \label{defx}
x := \alpha v^1 + \beta v^2.
\end{align}
It readily follows that $x \in \{ 0, \alpha, \beta \}^n.$
By using the substitution \eqref{defx}, we rewrite the objective as follows
\[
x^T Qx = \alpha^2 (v^1)^T Q v^1 + 2 \alpha \beta  (v^1)^T Q v^2  +    \beta^2 (v^2)^T Q v^2.
\]
Let us now introduce the following binary vector
$y^T = \begin{bmatrix}     (v^1)^T,     (v^2)^T    \end{bmatrix} \in \{ 0,1\}^{2n}$
and matrix
\begin{align}\label{barQ}
\overline{Q} =
\begin{bmatrix}
  \alpha^2 Q & \alpha \beta Q \\
  \alpha \beta Q & \beta^2 Q
\end{bmatrix} \in  {\mathcal S}^{2n}(\mathbb{R}),
\end{align}
from where it follows
$x^T Qx = y^T \overline{Q} y.$
Similarly, we can rewrite the constraints
\[
Ax = A(\alpha v^1 + \beta v^2) =
\begin{bmatrix}
    \alpha A, \beta A
\end{bmatrix} y
=
(\begin{bmatrix}
    \alpha,  \beta
\end{bmatrix} \otimes A) y
= b.
\]
The IQP \eqref{QP_with_3_values} can be now reformulated as a binary QP, and then by  exploiting Proposition~\ref{prop:Rank1} further reformulated into the following BSDP, see \cite{DeMeijerSotirov2024}:
\begin{equation} \label{BQP_with_3_valuesBin}
   \begin{array}{cl}
      \min & \langle \overline{Q}, Y \rangle  \\[1ex]
      \st    & (\begin{bmatrix}
    \alpha , \beta
     \end{bmatrix} \otimes A) y = b \\[1ex]
   &  \begin{bmatrix}
    I_n, I_n
     \end{bmatrix} y \leq e \\ [1ex]
     & \begin{bmatrix} 1 & y^T \\ y & Y \end{bmatrix} \succeq 0, ~~ \diag (Y) = y \\[2.5ex]
            & y \in \{ 0,1 \}^{2n}.
   \end{array}
\end{equation}
The binarity of the matrix variable $Y$ follows from the integrality of $y$ and the positive semidefiniteness constraint.

\medskip

\noindent {\bf Ternary  reformulation for not-equally spaced points.}
We may assume w.l.o.g.~that $\gamma=0$ and $0<\alpha<\beta$, $\beta \neq 2\alpha$.
The quadratic mapping $\varphi : \{-1,0,1\}^n  \to \{0, \alpha, \beta\}^n$ where
$\varphi(t) = \alpha e+\frac{\beta}{2}t
    +\left(\frac{\beta}{2}-\alpha\right)t\circ t$, maps $\{-1, 0,1 \}^n$ onto $\{0, \alpha, \beta\}^n.$
For $t_i\in\{-1,0,1\}$, define $z_i:=t_i^2=|t_i|$,
$w^T :=\begin{bmatrix}t^T, z^T\end{bmatrix}$ and
\begin{align}\label{matrixB}
    B:=\begin{bmatrix}
        \frac{\beta}{2}I_n,
        \left(\frac{\beta}{2}-\alpha\right)I_n
    \end{bmatrix}\in \mathbb{R}^{n\times 2n}.
\end{align}
Then for any $x\in \{0,\alpha,\beta\}^n$ we have $x=\alpha e+Bw$, and
the objective in \eqref{QP_with_3_values} can be rewritten as follows:
\begin{align}\label{ternaryObj}
x^TQx =\alpha^2e^TQe+2\alpha e^TQBw+w^TB^TQBw,
\end{align}
while  the equality constraint $Ax=b$ can be written as
\begin{align}\label{ternaryCostrEq}
ABw=b-\alpha Ae.
\end{align}
Thus, one  can reformulate the IQP \eqref{QP_with_3_values} by using \eqref{ternaryObj}--\eqref{ternaryCostrEq}  and
$z=|t|$ as an quadratic optimization problem with variables $t\in \{-1,0,1\}^n$, $z\in \{0,1\}^n.$

To reformulate the resulting problem as an ISDP, we proceed as follows.
 We introduce the lifted matrix
\begin{align}\label{blockM}
M =     \begin{bmatrix}
        1 & t^T & z^T\\
        t & T   & R\\
        z & R^T& Z
    \end{bmatrix} \succeq 0,
\end{align}
to which we impose the following constraints
\begin{align}\label{ternaryRank1}
    \diag(T)=z,\quad
    \diag(Z)=z,\quad
    z=|t|,\quad
    t\in\{-1,0,1\}^n.
\end{align}
These constraints together with the Schur complement lemma applied to $M$  imply that the rank of $M$ is one.
\begin{lemma}\label{rankOneTernary}
Let $M \in {\mathcal S}^{2n+1}$ be of the block form \eqref{blockM}   and satisfies the constraints~\eqref{ternaryRank1}. Then $\rank(M)=1.$
\end{lemma}
\begin{proof}
 It follows from the Schur complement applied to $M$ that
 $$
 \begin{bmatrix}
         T   & R\\
         R^T & Z
 \end{bmatrix}  -  \begin{bmatrix} t \\z\end{bmatrix}\begin{bmatrix} t \\z\end{bmatrix}^T   \succeq 0.
 $$
Then, it follows from $\diag \left (
\begin{bmatrix}
         T   & R\\
         R^T & Z
 \end{bmatrix} \right ) -
 \diag \left (  \begin{bmatrix}
         tt^T   & tz^T\\
         zt^T & zz^T
 \end{bmatrix} \right ) = \begin{bmatrix} z\\z \end{bmatrix} -  \begin{bmatrix} z\\z \end{bmatrix}=0$
 that
 $\begin{bmatrix}
         T   & R\\
         R^T & Z
 \end{bmatrix} = \begin{bmatrix}
t\\ z
\end{bmatrix}
\begin{bmatrix}
t\\ z
\end{bmatrix}^T $. Hence $\rank(M)=1.$
\end{proof}
The previous discussion results in the following
ISDP reformulation of \eqref{QP_with_3_values}:
\begin{equation} \label{ternary-ISDP}
\begin{array}{cl}
\displaystyle\min
    & \langle\widehat Q,M\rangle\\[1ex]
\st
    & \displaystyle
      AB\begin{bmatrix}t\\z\end{bmatrix}=b-\alpha Ae\\[3ex]
    & \displaystyle
      M=      \begin{bmatrix}
          1 & t^T & z^T\\
          t & T   & R\\
          z & R^T & Z
      \end{bmatrix}\succeq0 \\[3.5ex]
    & \diag(T)=z,\quad \diag(Z)=z \\[1ex]
    & z=|t|, \quad  t\in\{-1,0,1\}^n, \quad z \in\{0,1\}^n,
\end{array}
\end{equation}
where
\begin{align}\label{widehatQ}
\widehat Q:=     \begin{bmatrix}
        \alpha^2e^TQe & \alpha e^TQB\\
        \alpha B^TQe  & B^TQB
    \end{bmatrix} \in  {\mathcal S}^{2n+1}(\mathbb{R}).
    \end{align}
Note that the matrix variables in \eqref{BQP_with_3_valuesBin} and \eqref{ternary-ISDP} both have order $2n+1.$
Integrality of $T,R,Z$ follows from the integrality of $t,z$ and rank-one property of $M$.
Note also that the ISDP~\eqref{ternary-ISDP} contains  a nonconvex constraint $ z=|t|$.
The idea of reformulating ternary QPs as ternary SDPs was introduced in~\cite{deMeijer2026BeyondBinarity}.

Observe that when $\beta = 2\alpha$, the ISDP~\eqref{ternary-ISDP} reduces to a (rescaled) ternary SDP.

\medskip\medskip\medskip

\noindent
{\bf Ternary reformulation for equally spaced points.}
Suppose that  $x\in \{-\alpha, 0, \alpha\}^n$ for some  $\alpha \in \mathbb{Z}$, $\alpha > 0$   in the IQP \eqref{QP_with_3_values}.
Then one can apply Lemma~\ref{lemma:alphaSminalpha0-1}
to reformulate the IQP~\eqref{QP_with_3_values} as an ISDP whose matrix variable has entries in the set $\{ -\alpha^2, 0, \alpha^2\}.$
In particular the IQP~\eqref{QP_with_3_values} can be equivalently formulated as:
\begin{equation} \label{ISDP_with_3_values}
   \begin{array}{cl}
      \min & \langle Q, X \rangle \\[1ex]
      \st    & A x = b \\[1ex]
              & \begin{bmatrix} 1 & x^T \\ x & X \end{bmatrix} \succeq 0, ~~ \diag (X) = \alpha |x| \\[2ex]
            & X \in \{-\alpha^2, 0, \alpha^2 \}^{n \times n}.
   \end{array}
\end{equation}
Note that the ternary SDP \eqref{ISDP_with_3_values} involves a matrix variable of order $n+1$, while the binary SDP~\eqref{BQP_with_3_valuesBin} requires a matrix variable of order $2n+1.$ Here we also allow additional nonlinear constraint $\diag (X) = \alpha |x|$.
Solving ternary SDPs can be efficiently done for matrix variables of order 150, see~\cite{deMeijer2026BeyondBinarity}.
\begin{remark}
Observe that every equally spaced three-point set
$\{\delta, (\delta+\beta)/2, \beta\}$ can be mapped bijectively onto $\{-\alpha, 0, \alpha\}$ for any $\alpha >0$, and thus reformulated as the ISDP~\eqref{ISDP_with_3_values}.
\end{remark}

\subsection{Comparison of SDP relaxations for not-equally spaced points}
\label{sect:compareSDPThreeGeneral}
We compare SDP relaxations obtained from  \eqref{BQP_with_3_valuesBin} and \eqref{ternary-ISDP} for
not-equally spaced points. That is, we assume that $\gamma=0$ and $0<\alpha<\beta$, $\beta \neq 2\alpha$. We prove that these two SDP relaxations are equivalent.

One can derive a basic SDP relaxation from the BSDP~\eqref{BQP_with_3_valuesBin} by relaxing binary constraints and adding  sparsity  on the matrix variable.
The sparsity pattern follows from the fact that $v^1_i + v^2_i \leq 1$ for all $i\in [n],$ see \eqref{defx}. Those constraints are automatically satisfied in \eqref{BQP_with_3_valuesBin}, but not in its relaxation. In particular, the BSDP \eqref{BQP_with_3_valuesBin} gives a rise for the following SDP relaxation:
\begin{equation} \label{SDPrelax-BQP_with_3_valuesBin}
   \begin{array}{cl}
      \min & \langle \overline{Q}, Y \rangle  \\[1ex]
      \st    & \left (\begin{bmatrix}
    \alpha , \beta
     \end{bmatrix} \otimes A \right ) y = b \\[1ex]
   &  \begin{bmatrix}
    I_n, I_n
     \end{bmatrix} y \leq e \\ [1ex]
      & Y_{i,n+i} = 0, \quad Y_{n+i,i} = 0 \qquad i\in [n] \\[1ex]
     & \widetilde{Y} = \begin{bmatrix} 1 & y^T \\ y & Y \end{bmatrix} \succeq 0, ~~ \diag (Y) = y, ~~\\[2.ex]
     & y\in\mathbb R^{2n}, ~~Y\in {\mathcal S}^{2n},
   \end{array}
\end{equation}
where  $\overline{Q}$ is given in \eqref{barQ}.
We will further use the following decomposition of $Y\in {\mathcal S}^{2n}$ and $y\in \mathbb{R}^{2n}$:
\begin{align}  \label{decomposition}
Y = \begin{bmatrix} Y^{11} & Y^{12} \\ Y^{21} & Y^{22} \end{bmatrix},
\qquad y= \begin{bmatrix} y^1\\ y^2 \end{bmatrix},
\end{align}
where $ Y^{12} = ( Y^{21})^T \in \mathbb{R}^{n\times n},$ $Y^{11}, Y^{22} \in \mathcal{S}^n$, and  $y^1, y^2\in \mathbb{R}^n$.

After relaxing the integrality conditions in the ISDP \eqref{ternary-ISDP}, and adding a constraint on the diagonal of $R$, we obtain the following SDP relaxation for \eqref{QP_with_3_values}:
\begin{equation}
\begin{array}{cl}  \label{SDPrelax-BQP_with_3_valuesTer}
\min
    & \langle\widehat Q,M\rangle\\[1ex]
\st
    &       AB\begin{bmatrix}t\\z\end{bmatrix}=b-\alpha Ae\\[3ex]
    &       M=
      \begin{bmatrix}
          1 & t^T & z^T\\
          t & T   & R\\
          z & R^T & Z
      \end{bmatrix}\succeq0 \\[3.5ex]
    & \diag(T) = z, \quad \diag(Z) = z, \quad \diag(R)=t \\[1ex]
    &  t,z\in\mathbb R^n,\quad
      T,Z\in\mathcal S^n,\quad
      R\in\mathbb R^{n\times n},
\end{array}
\end{equation}
where $\widehat Q$ is given in \eqref{widehatQ} and $B$ is defined in \eqref{matrixB}.
The constraint $ \diag(R)=t $ follows from the fact that  $R=t z^T$ and $z=|t|$, $t\in \{-1, 0,1\}^n,$ in \eqref{ternary-ISDP}.
It follows from $M\succeq 0$ and $\diag(Z)=z$ that $\diag(T)=\diag(Z)\leq e_n.$
Observe that the diagonal constraints on $M$ and $M\succeq 0$ imply that
$z \geq t,$ $z\geq -t .$
The following result relates the above two SDP relaxations.

\begin{theorem} \label{thm:SDPequivalence-Tri}
Let $\alpha, \beta \in \mathbb{Z}$, $0<\alpha<\beta$, $\beta \neq 2\alpha$.
Then, the SDP relaxations~\eqref{SDPrelax-BQP_with_3_valuesBin}
 and~\eqref{SDPrelax-BQP_with_3_valuesTer}  are equivalent.
\end{theorem}
\begin{proof}
Let $(y,Y)$ be feasible for \eqref{SDPrelax-BQP_with_3_valuesBin}.
Assume that $(y,Y)$ is partitioned as in~\eqref{decomposition}.
Define  $M := P \widetilde{Y} P^T$ where
\begin{align}\label{matrixP}
P:= \begin{bmatrix}
1 & 0 & 0\\
-e & I_n & 2I_n\\
e &-I_n & 0
\end{bmatrix}.
\end{align}
Assume that $M$ is partitioned as in \eqref{blockM}.
We show that $M$ is feasible for~\eqref{SDPrelax-BQP_with_3_valuesTer}.
Since $\widetilde Y\succeq0$, it follows that $M\succeq0$.

Since the  first row of $M$ equals  $\begin{bmatrix} 1, ~ (-e+y^1+2y^2)^T, ~(e-y^1)^T \end{bmatrix},$
we set $t =-e+y^1+2y^2$,  $z =e-y^1$.
Let us  verify the diagonal constraints of $M$. The blocks of $M$ obtained by multiplying out are
\begin{align} \label{blockB}
T & =  J_n -e \left( y^1 + 2 y^2\right )^T - \left ( y^1 + 2y^2 \right ) e^T + Y^{11} + 2Y^{12} + 2 Y^{21} + 4Y^{22} \\
Z & = J_n -e(y^1)^T-y^1 e^T+Y^{11} \label{blockZ}\\
R & = -J_n +e(y^1)^T+ \left ( y^1 + 2y^2 \right ) e^T -Y^{11}-2Y^{21}. \label{blockR}
\end{align}
It follows from \eqref{blockB}--\eqref{blockR}, after  exploiting
$\diag(Y)=y$ and $\diag(Y^{12})=\diag(Y^{21})=0,$ that
$\diag(T) =z,$ $\diag(Z)=z,$ $\diag(R)=t.$

To verify equality constraints, we proceed  as follows:
\begin{align*}
A \left ( \alpha e+B\begin{bmatrix}t\\z\end{bmatrix} \right )
& = A \left(\alpha e +\frac{\beta}{2}(-e+y^1 + 2y^2)
+\left(\frac{\beta}{2}-\alpha\right)(e-y^1)\right)
=A \left ( \alpha y^1 + \beta y^2 \right)=b.
\end{align*}
Thus, $M$ is feasible for~\eqref{SDPrelax-BQP_with_3_valuesTer}.
\medskip

Conversely,  let $M$  be feasible for \eqref{SDPrelax-BQP_with_3_valuesTer}.
Define   $K:= P^{-1}$, see \eqref{matrixP}, i.e.,
$$
K:=
\begin{bmatrix}
1 & 0& 0\\
e & 0 &-I_n\\
0 &\frac{1}{2} I_n &\frac{1}{2} I_n
\end{bmatrix}.
$$
Define
$$
\widetilde{Y} := K M K^T =
\begin{bmatrix}
1 & (y^1)^T & (y^2)^T \\
y^1 & Y^{11} & Y^{12} \\
y^2 & Y^{21} & Y^{22}
\end{bmatrix}.
$$
We show below that $\widetilde{Y}$ is feasible for \eqref{SDPrelax-BQP_with_3_valuesBin}.
Since $M\succeq0$, it follows that $\widetilde{Y} \succeq0$.
A direct calculation gives
$y^1 = e-z$, $y^2=\frac{1}{2}(t+z)$, and
\begin{align*}
Y^{11}  = J_n-e z^T - z e^T+Z, \,  
Y^{12} =\frac{1}{2} \left( e (t + z)^T - R^T-Z \right),  \,
Y^{22}  = \frac{1}{4} \left( T+R+R^T+Z \right ).
\end{align*}
Using
$\diag(T)=\diag(Z)=z,$ $\diag(R)=t$, we obtain
$\diag(Y^{11})=e-z=y^1,$
$\diag(Y^{22}) =\frac{1}{2}(t+z)=y^2.$ Thus $\diag(Y)=y.$

Also $\diag(Y^{12})=\diag(Y^{21}) =\frac{1}{2}(t+z-t-z)=0.$ Furthermore,
$$
y^1 + y^2 = e-z+\frac{1}{2}(t+z) = e + \frac{t-z}{2} \le e,
$$
since $t-z \leq 0.$
We verify the remaining equality constraints by exploiting
$$
\left( \begin{bmatrix}  \alpha, \beta  \end{bmatrix} \otimes A \right ) y
= \alpha Ay^1 + \beta Ay^2 = \alpha A(e-z)+\frac{\beta}{2} A (t+z)
= \alpha Ae + AB\begin{bmatrix}t\\z\end{bmatrix} =b.
$$
It remains to verify that the two objectives are equal. Define
\begin{align} \label{CD}
C:=\begin{bmatrix}\alpha I_n, ~\beta I_n\end{bmatrix},
\qquad
D := \begin{bmatrix}\alpha e, ~B \end{bmatrix}.
\end{align}
Note that $C$ and $D$ are defined such that
$\widehat{Q} = D^T Q D$ and $\overline{Q} = C^T Q C.$

Since $M =P \widetilde{Y} P^T$ and  $\widetilde{Y} := K M K^T$, we have
\begin{align*}
D M D^T
 = D \left ( P\widetilde{Y} P^T \right ) D^T =
\begin{bmatrix}0, ~ C\end{bmatrix}
\widetilde{Y}
\begin{bmatrix}0, ~C \end{bmatrix}^T
 = CYC^T,
\end{align*}
where we exploited $DP =  \begin{bmatrix}0,~  C \end{bmatrix}.$
Finally,
$\langle \widehat{Q},M\rangle
=\langle D^T Q D,M\rangle
=\langle Q,D M D^T\rangle
=\langle Q,CYC^T\rangle
=\langle C^TQC,Y\rangle
=\langle \overline{Q},Y\rangle.$
\end{proof}
Observe that the SDP relaxation~\eqref{SDPrelax-BQP_with_3_valuesBin} has $m+4n$ equality and $n$ inequality constraints, and the SDP relaxation~\eqref{SDPrelax-BQP_with_3_valuesTer} has $m+3n$ equality  constraints, while both have matrix variables of order $2n+1.$
Thus, the SDP relaxation~\eqref{SDPrelax-BQP_with_3_valuesTer} is more favorable for computing bounds.

\subsection{Comparison of SDP relaxations for equally-spaced points}
\label{sect:SDP relax of equally spaced}
We compare SDP relaxations obtained from ISDPs \eqref{BQP_with_3_valuesBin} and
\eqref{ISDP_with_3_values} for the case    $x\in \{-\alpha, 0, \alpha\}^n$ where
$\alpha \in \mathbb{Z}$, $\alpha > 0$.
The results in this section show that a ternary SDP relaxation derived from  \eqref{ISDP_with_3_values} is more favorable than a  binary SDP relaxation of \eqref{BQP_with_3_valuesBin}.
Namely, the former relaxation has a matrix variable of smaller order and fewer constraints than the latter, while providing the same bounds.

We consider the following SDP relaxation that is obtained from the ternary SDP \eqref{ISDP_with_3_values}:
\begin{equation} \label{SDPrelax - ISDP_with_3_values}
   \begin{array}{cl}
      \min & \langle Q, X \rangle \\[1ex]
      \st    & A x = b \\[1ex]
      & \diag (X) \geq \alpha x , \quad   \diag (X) \geq - \alpha x \\[1ex]
              &  \diag (X) \leq \alpha^2 e, \quad \begin{bmatrix} 1 & x^T \\ x & X \end{bmatrix} \succeq 0,
   \end{array}
\end{equation}
where $X\in  {\mathcal S}^n(\mathbb{R}).$
In what follows, we show that the SDP relaxations \eqref{SDPrelax-BQP_with_3_valuesBin} and  \eqref{SDPrelax - ISDP_with_3_values} are equivalent. We first prove the easier inclusion.
\begin{proposition} \label{SDPequivalence-first part}
Let  $\alpha \in \mathbb{Z}$, $\alpha > 0$, and $\beta = -\alpha$.
Then, every feasible solution of the binary SDP relaxation~\eqref{SDPrelax-BQP_with_3_valuesBin} projects to a feasible solution of the ternary SDP relaxation~\eqref{SDPrelax - ISDP_with_3_values} with the same objective value.
\end{proposition}
\begin{proof}
Let $(y,Y)$ be feasible for the  SDP relaxation~\eqref{SDPrelax-BQP_with_3_valuesBin}.
We define $W:=\alpha[I_n,-I_n].$
Let
$$x:=W y=\alpha(y^1-y^2), \quad X:= WY W^T.$$
In what follows, we verify that $(x,X)$ is feasible for the SDP relaxation~\eqref{SDPrelax - ISDP_with_3_values}.
First, we show that  the linear equality constraints hold:
\[
Ax = A(W y) =
\begin{bmatrix}
\alpha A, ~ -\alpha A
\end{bmatrix} y =
(\begin{bmatrix} \alpha, -\alpha \end{bmatrix}\otimes A)y=
b.
\]
Next, we verify the diagonal inequalities by exploiting the matrix block-decomposition~\eqref{decomposition}:
\[
\diag(X) = \diag(W Y W^T)
= \alpha^2 \diag( Y^{11} + Y^{22} - Y^{12} - Y^{21})=
\alpha^2 (y^1 + y^2) \leq \alpha^2 e.
\]
Here, we exploit  the fact that $\diag(Y^{12} + Y^{21}) =0$ and
$\begin{bmatrix}  I_n, I_n   \end{bmatrix} y \leq e$.
Further, since $y=\diag(Y) \geq 0$ we have
$$
X_{ii} = \alpha^2 (y_i^1 + y_i^2) \geq \alpha^2 | y_i^1 - y_i^2| = \alpha |x_i|, \qquad \forall i\in [n],
$$
which implies the remaining diagonal inequalities.

Let us verify the PSD constraint. Define
$R := \begin{bmatrix} 1 & 0\\ 0 & W \end{bmatrix}.$
Then
$$
R \begin{bmatrix}
1 & y^T\\
y & Y
\end{bmatrix}
R^T =
\begin{bmatrix}
1 & x^T\\
x & X
\end{bmatrix} \succeq 0.
$$
It remains to verify that the two objectives are equal. Indeed
$\langle Q, X \rangle = \langle Q, W Y W^T \rangle  =
\langle W^T Q W, Y \rangle = \langle \overline{Q}, Y \rangle,$
where we used that  $\overline{Q}= W^T Q W$, see \eqref{barQ}.
\end{proof}
We show now the reverse statement.
\begin{proposition} \label{SDPequivalence-second part}
Let   $\alpha \in \mathbb{Z}$, $\alpha > 0$, and $\beta=-\alpha$.
Then, every feasible solution of the  SDP relaxation~\eqref{SDPrelax - ISDP_with_3_values} can be lifted to a feasible solution of the  SDP relaxation~\eqref{SDPrelax-BQP_with_3_valuesBin} with the same objective value.
\end{proposition}
\begin{proof}
 Let $(x,X)$ be feasible for the SDP relaxation~\eqref{SDPrelax - ISDP_with_3_values}.
Let  $z :=\frac{1}{\alpha}x,$ $Z:=\frac{1}{\alpha^2}X.$
Then, it readily follows that $\diag(Z)\ge z,$ $\diag(Z)\ge -z,$ $\diag(Z)\le e.$ Moreover, it follows from $X- xx^T\succeq 0$ that
$$ \begin{bmatrix} 1 & z^\top\\ z & Z \end{bmatrix}\succeq 0.$$
Define
$
d:=\diag(Z),$
$y^1 := \frac{1}{2}(d+z),$
$y^2 := \frac{1}{2}(d - z).$
Then
$$
y^1, y^2 \ge 0,\qquad
y^1 + y^2 = d \le e, \qquad
y^1 - y^2 = z.
$$
Hence  $x=\alpha z = W y$ where
$W:=\alpha [I_n,-I_n]$ and
$y^T= \begin{bmatrix} (y^1)^T, (y^2)^T \end{bmatrix}$.
We verify the equality constraint:
$
([\alpha,-\alpha]\otimes A)y
=
A\alpha(y^1 - y^2)
= Ax=b.
$

To obtain a lifted matrix we proceed as follows.
Let
$
S := Z-zz^T \succeq 0,
$
and $h^1,\ldots,h^n \in \R^{2n}$ be vectors satisfying
$$
S_{ij}=(h^i)^T h^j \qquad \forall i,j \in [n].
$$
Here we assume that,  for each $i\in [n]$, the first $n$ components of $h^i$  are obtained, for instance, from an eigendecomposition, and the remaining components are zero.

In the sequel, we introduce vectors $g_1^{i}$, $g_2^{i}$ for $i\in [n]$
that satisfy  $h^i = g_1^{i} - g_2^{i}.$
For each $i \in [n]$, let $s_i := d_i-z_i^2 = ||h^i ||^2$.
If $s_i=0$, set $g_1^{i}=g_2^{i}=0$.
Otherwise, let $u^i$ for $i\in [n]$ denotes the vector in $\R^{2n}$  obtained by embedding the standard unit vector $u_i$ of $\R^n$ into the last $n$ coordinates. Thus $u^i$ is orthogonal to $h^i.$
Define
$$
g_1^{i} := a_i h^i + t_i u^i,\quad
g_2^{i}:= (a_i-1)h^i + t_i u^i,
$$
where
$$
a_i =\frac{1}{2}\left(1+\frac{z_i(1-d_i)}{s_i}\right),
\quad
t_i^2 = \frac{(d^2_i-z^2_i)(1-d_i)} {4s_i}.
$$
For so defined vectors, we have that  $g_1^{i}-g_2^{i}=h^i,$
and
\begin{align}\label{d1g2}
\|g_1^{i}\|^2 = y_i^{1}(1-y_i^{1}),\quad
\|g_2^{i}\|^2 =y_i^{2}(1-y_i^{2}).
\end{align}
Let $G:=\begin{bmatrix} g_1^{1}, g_1^{2},\ldots, g_1^{n}, g_2^{1}, g_2^{2},\ldots, g_2^{n}\end{bmatrix}$
and
\begin{align}\label{defY}
Y= yy^T + G^T G.
\end{align}
Then
$Y-yy^T =  G^T G \succeq 0,$ which is equivalent to
$\begin{bmatrix}
1 & y^T \\
y &  Y
\end{bmatrix} \succeq 0.$

It follows from \eqref{d1g2} and \eqref{defY} that
$Y_{ii} =(y_i^{1})^2+\|g_1^{i}\|^2=y_i^1,$
$Y_{n+i,n+i}=  (y_i^{2})^2+ \|g_2^{i}\|^2 = y_i^2.$
Thus, $\diag(Y)=y.$
It follows from \eqref{d1g2} and \eqref{defY} that
$$
Y_{i,n+i} = y_i^1 y_i^2 + (g_1^{i})^T g_2^{i} = \frac{d_i^2-z_i^2}{4} - \frac{d_i^2-z_i^2}{4}=0,
$$
and similarly $Y_{n+i,i}=0.$
Thus, the sparsity constraints are also satisfied.

It remains to verify the objective.
We have
$$
W Y W^T = Wyy^T W^T  + W(Y-yy^\top) W^T = \alpha^2 zz^T + \alpha^2 S = \alpha^2 Z = X,
$$
where $W=\alpha [I_n,-I_n]$.
The second equality above follows from:
$$
\bigl(W(Y-yy^T)W^T\bigr)_{ij}
= \alpha^2 (g_1^i-g_2^i)^T(g_1^j-g_2^j)\\
= \alpha^2 (h^i)^T h^j\\
= \alpha^2 S_{ij}.
$$
Since $\overline{Q} = W^TQW$, it follows that
$\langle \overline{Q},Y\rangle = \langle W^TQW,Y\rangle = \langle Q,W Y W^T\rangle
= \langle Q,X\rangle.$
\end{proof}
The following results follows from Proposition~\ref{SDPequivalence-first part} and  Proposition~\ref{SDPequivalence-second part}.
\begin{theorem}
Let   $\alpha \in \mathbb{Z}$, $\alpha > 0$, and $\beta =-\alpha$.
Then, the  SDP relaxation~\eqref{SDPrelax - ISDP_with_3_values}
is equivalent to the SDP relaxation~\eqref{SDPrelax-BQP_with_3_valuesBin}.
\end{theorem}
We remark that the  binary SDP relaxation~\eqref{SDPrelax-BQP_with_3_valuesBin} without the sparsity constraints is dominated by the
ternary SDP relaxation~\eqref{SDPrelax - ISDP_with_3_values}.

\section{IQPs with a finite number of integer values as ISDPs}
\label{sect:generalISDP}
We consider the following QP with a constraint that the variable can only take integer values in a finite range $[L,U]$ of values. Since we can always apply a simple coordinate shift, it is no loss of generality to assume that $L = 0$. So to be more specific, we consider the following IQP:
\begin{equation} \label{pure_integer_QP}
   \begin{array}{cl}
      \min & x^T Q x \\
      \st    & Ax = b \\
              & x \in \Z^n \cap [ 0, U ]^n,
   \end{array}
\end{equation}
where $Q \in {\mathcal S}^n(\mathbb{R})$, $A \in \R^{ m \times n }$, $b \in \R^m$ and $U \in \Z_+$.

We first drive three equivalent integer SDP formulations of the IQP \eqref{pure_integer_QP},
and then in Section~\ref{sect:comparision of general formularins} compare the corresponding SDP relaxations.
The results show that the SDP relaxation of the ISDP formulation with the largest matrix variable, that is the one obtained by one-hot decoding, dominates other two SDP relaxations.
However, that relaxation may become too computationally expensive to solve as $n$ or the number of distinct integer values increases.

\medskip
\noindent {\bf One-hot decoding.}
Let us define
\begin{align}\label{defu}
\mathfrak{u} := \begin{bmatrix} 1, 2, \ldots, U \end{bmatrix}^T \in \Z^U.
\end{align}
In order to reformulate \eqref{pure_integer_QP} as a BQP, we introduce $nU$ binary variables as follows: 
\begin{align}  \label{vij}
v^i_j := \begin{cases}
                0 & \text{ if } x_i \neq j \\
                1 & \text{ if } x_i = j,
            \end{cases}
\end{align}
where $i \in [n]$ and $j \in [U]$.
This gives rise to a matrix $V \in \R^{n \times U}$ defined by $V_{ij} := v^i_j$. It is then easy to verify that $x = V \mathfrak{u}$. To ensure this representation of $x$ is unique, we impose the one-hot constraint
\begin{equation} \label{eq:Ve leq e}
    Ve \leq e.
\end{equation}
Applying the $\vect$ operator to the matrix $V$, it is not difficult to see that the one-hot decoding $x = V \mathfrak{u}$ can be rewritten as
\begin{equation}\label{binEncoding_xi}
    x = (I_n \otimes \mathfrak{u}^T) \vect (V^T).
\end{equation}
With this, we can rewrite the objective and constraints of \eqref{pure_integer_QP} and  arrive at
\begin{equation} \label{binary_formulation_of_pure_integer_QP}
   \begin{array}{cl}
      \min & \vect (V^T)^T (I_n \otimes \mathfrak{u}^T)^T Q (I_n \otimes \mathfrak{u}^T) \vect (V^T) \\[1ex]
      \st   & A (I_n \otimes \mathfrak{u}^T) \vect (V^T) = b \\[1ex]
            & (I_n \otimes e^T) \vect (V^T) \leq e \\[1ex]
            & \vect (V^T) \in \{ 0,1 \}^{nU}.
   \end{array}
\end{equation}
Now, we define $y := \vect (V^T)$, and reformulate the binary QP \eqref{binary_formulation_of_pure_integer_QP}  by using  Proposition~\ref{prop:Rank1}, see also~\cite{DeMeijerSotirov2024}:
\begin{equation} \label{SDP_formulation_of_pure_integer_QP}
   \begin{array}{cl}
      \min & \left\langle (I_n \otimes \mathfrak{u}) Q (I_n \otimes \mathfrak{u}^T), Y \right\rangle \\[1ex]
      \st   & A (I_n \otimes \mathfrak{u}^T) y = b \\[1ex]
            & (I_n \otimes e^T) y \leq e \\[1ex]
            & \begin{bmatrix} 1 & y^T \\ y & Y \end{bmatrix} \succeq 0, ~~ \diag (Y) = y \\[2ex]
            & y \in \{ 0,1 \}^{nU}.
   \end{array}
\end{equation}
This ISDP has a matrix variable of order $nU+1.$
We show in the following how to derive ISDP reformulations of \eqref{pure_integer_QP} that have  matrix  variables of order smaller  than $nU+1$.

\medskip
\noindent {\bf Binary expansion.}
We  first rewrite \eqref{pure_integer_QP} as an equivalent binary  quadratic program.
Each \(x_i\in \mathbb Z\cap[0,U]\), \(i\in[n]\), can be uniquely represented as
\begin{align}\label{sumx2}
x_i=\sum_{j=1}^l 2^{l-j}b_j^i,
\qquad
b_j^i\in\{0,1\},\quad  j \in[l],
\end{align}
where
$l =\left\lceil \log_2(U+1)\right\rceil.$
Define $B \in\mathbb \{0,1\}^{n\times l}$ by $B_{ij}:= b_j^i\in\{0,1\},$ for $i \in[n],\ j\in[l],$
and
\begin{align}\label{binaryv}
\mathfrak{w}:= \begin{bmatrix} 2^{l-1}, 2^{l-2}, \cdots,  1 \end{bmatrix}^T. \end{align}
Then $x= B\mathfrak{w},$ which can equivalently be written as
$x= (I_n\otimes\mathfrak{w}^T)\vect(B^T).$

Using this equality, we rewrite the objective and the constraints of \eqref{pure_integer_QP}, which result~in:
$$
\begin{array}{cl}
\min &
\vect(B^T)^T (I_n\otimes\mathfrak{w})Q
(I_n\otimes\mathfrak{w}^T) \vect(B^T) \\[1ex]
\st &
A(I_n\otimes\mathfrak{w}^T)\vect(B^T)=b
\\[1ex]
& (I_n\otimes\mathfrak{w}^T)\vect(B^T)\leq Ue
\\[1ex]
& \vect(B^T)\in\{0,1\}^{nl}.
\end{array}
$$
Here we add the upper bound on $x=(I_n\otimes\mathfrak{w}^T)\vect(B^T)$ to exclude values represented by the binary decomposition that are larger than $U.$

Let $y:= \vect(B^T)$. Now, by using Proposition~\ref{prop:Rank1}, we obtain the following binary SDP reformulation of the IQP \eqref{pure_integer_QP}:
\begin{align} \label{SDPexactBinaryEnc}
\begin{array}{cl}
\min &
\left\langle
(I_n \otimes \mathfrak{w}) Q (I_n \otimes \mathfrak{w}^T),\, Y
\right\rangle
\\[1ex]
\st & A (I_n \otimes \mathfrak{w}^T) y = b
\\[1ex]
& (I_n \otimes \mathfrak{w}^T) y \leq Ue
\\[1ex]
& \begin{bmatrix}
1 & y^T \\
y & Y
\end{bmatrix} \succeq 0,
\quad \diag(Y) =y \\[2ex]
& y\in\{0,1\}^{nl}.
\end{array}
\end{align}
Since $l \leq U$, the BSDP \eqref{SDPexactBinaryEnc} has a smaller matrix variable than
the BSDP \eqref{SDP_formulation_of_pure_integer_QP}.
We can further reduce the number of variables by using a ternary number representation.

\medskip
\noindent {\bf Ternary expansion.}
We proceed similarly as for the binary expansion.
First observe that each $x_i \in \Z \cap [ 0, U ]$ for $i\in [n]$ can be uniquely represented as:
\begin{align}    \label{sumx3}
x_i = \sum_{j=1}^k 3^{k-j} ( t^i_j +1), \qquad t^i_j \in \{-1,0,1\}, \,\, j\in [k],
\end{align}
where $k= \left \lceil ~\log_3(U+1)~ \right \rceil.$
Define $T \in \{-1,0,1\}^{n\times k}$ such that $T_{ij} := t^i_j \in
\{-1,0,1\}$ for $i\in [n]$, $j\in [k]$, and
\begin{align}\label{mathfrakv}
\mathfrak{v}:= \begin{bmatrix} 3^{k-1}, 3^{k-2}, \cdots, 1 \end{bmatrix}^T.
\end{align}
Then, by exploiting $\sum_{j=1}^k 3^{k-j} = (3^k -1)/2$, it follows
$x = T\mathfrak{v} + (3^k -1)e/2$, which can be rewritten as
\begin{align*}
    x = (I_n \otimes \mathfrak{v}^T) \vect (T^T) + \frac{(3^k - 1)}{2}e.
\end{align*}
By proceeding as in the case of the binary expansion, we obtain the following ternary SDP:
\begin{equation} \label{ternarySDP_formulation_of_pure_integer_QP}
   \begin{array}{cl}
      \min &   \left\langle (I_n \otimes \mathfrak{v}) Q (I_n \otimes \mathfrak{v}^T), Y \right\rangle
       +        (3^k-1)e^T  Q (I_n \otimes \mathfrak{v}^T) y + \frac{(3^k-1)^2}{4}e^TQe \\[1ex]
      \st   & A (I_n \otimes \mathfrak{v}^T) y = b - \frac{3^k - 1}{2} A e \\[1ex]
            & (I_n \otimes \mathfrak{v}^T) y  \leq \left (U-\frac{(3^k - 1)}{2} \right )e  \ \\[1ex]
            & \begin{bmatrix} 1 & y^T \\ y & Y \end{bmatrix} \succeq 0, ~~ \diag (Y) = |y|,  ~~ y \in \{-1, 0,1 \}^{nk}.
   \end{array}
\end{equation}
Note that $k$ can be much smaller than $U$, and therefore the matrix variable in \eqref{ternarySDP_formulation_of_pure_integer_QP} can be much smaller than that in \eqref{SDP_formulation_of_pure_integer_QP}.
The integrality of $Y$ follows from the integrality of $y$, see~\cite{deMeijer2026BeyondBinarity}.
 Since also  $k\leq l$, the TSDP \eqref{ternarySDP_formulation_of_pure_integer_QP} has the smallest matrix variable among all  three reformulations of \eqref{pure_integer_QP}.

\subsection{Comparison of relaxations}
\label{sect:comparision of general formularins}
Here we compare SDP relaxations that are obtained from the BSDP~\eqref{SDP_formulation_of_pure_integer_QP},
BSDP \eqref{SDPexactBinaryEnc}, and TSDP~\eqref{ternarySDP_formulation_of_pure_integer_QP}.
The results show that the SDP relaxation obtained from the one-hot decoding dominates the other two SDP relaxations, which are incomparable.

The BSDP \eqref{SDP_formulation_of_pure_integer_QP}
gives  rise for the following SDP relaxation:
\begin{align} \label{BSDPrelax}
   \begin{array}{cl}
      \min & \left\langle (I_n \otimes \mathfrak{u}) Q (I_n \otimes \mathfrak{u}^T), Y \right\rangle \\[1ex]
      \st   & A (I_n \otimes \mathfrak{u}^T) y = b \\[1ex]
            & (I_n \otimes e^T) y \leq e \\[1ex]
          & Y_{(i-1)U+r, (i-1)U+s } = 0\qquad i\in[n],  \quad r,s\in[U],\quad r\neq s \\[1ex]
          &  \widetilde{Y} = \begin{bmatrix} 1 & y^T \\ y & Y \end{bmatrix} \succeq 0, ~~ \diag (Y) = y, 
   \end{array}
\end{align}
where $Y \in {\mathcal S}^{nU}$ and $\mathfrak{u}$ defined in \eqref{defu}.
In what follows,  we use the decomposition of  $Y \in {\mathcal S}^{nU \times nU}$ into $n^2$ blocks of size $U\times U$, and
denote those blocks by $Y^{ij}$ for $i,j\in [n]$.
Similarly we use the decomposition $y^T = \begin{bmatrix} (y^1)^T, \cdots, (y^{n})^T \end{bmatrix}$ where $y^i \in \mathbb{R}^U$ for $i\in [n]$.  By exploiting the block decomposition, the sparsity pattern in \eqref{BSDPrelax} can be rewritten as
$ Y^{ii} = \Diag(y^i)$  for $i\in [n],$
Sparsity constraints follow from the fact that $v^i(v^i)^T$ is a diagonal matrix, see~\eqref{vij}.

BSDP \eqref{SDPexactBinaryEnc} gives rise to the following SDP relaxation:
\begin{align}\label{BinExpSDPrelax}
\begin{array}{cl}
\min & \left\langle (I_n \otimes \mathfrak{w}) Q (I_n \otimes \mathfrak{w}^T),\, X \right\rangle \\[1ex]
\st & A (I_n \otimes \mathfrak{w}^T)x = b \\[1ex]
& (I_n \otimes \mathfrak{w}^T) x \leq U e \\[1ex]
& \begin{bmatrix}
1 & x^T \\
x & X
\end{bmatrix} \succeq 0,
\quad \diag(X)= x,
\end{array}
\end{align}
where $Y \in {\mathcal S}^{nl}$ and $\mathfrak{w}$ defined in \eqref{binaryv}.
Next results show that the SDP relaxation \eqref{BinExpSDPrelax} is dominated by the SDP relaxation~\eqref{BSDPrelax}.
\begin{theorem}
Let $U \in \Z_+$ and  $l=\lceil\log_2(U+1)\rceil$.
 Then, the one-hot SDP relaxation~\eqref{BSDPrelax} dominates the binary SDP relaxation~\eqref{BinExpSDPrelax}.
\end{theorem}
\begin{proof}
Let ${\mathcal B}\in\{0,1\}^{l\times U}$ be the matrix whose $r$th column is the binary representation of $r$, for $r\in[U]$. Thus, $ \mathfrak{w}^T {\mathcal B} =\mathfrak{u}^T.$

Let $(y,Y)$ be a feasible point of~\eqref{BSDPrelax}.
Define $x:=Dy$ and $X:=DYD^T,$ where  $D:=I_n\otimes {\mathcal B}$.
Then
$$
\begin{bmatrix}
1 & x^T\\
x & X
\end{bmatrix}
= \begin{bmatrix}
1 & 0\\
0 & D
\end{bmatrix}
\begin{bmatrix}
1 & y^T\\
y & Y
\end{bmatrix}
\begin{bmatrix}
1 & 0\\
0 & D^T
\end{bmatrix}
\succeq0.
$$
From the zero pattern of $Y$ and  $\diag(Y)=y$, it follows
$$
X_{(i-1)l+j,(i-1)l+j} =\sum_{r=1}^U {\mathcal B}_{jr}^2 y_{(i-1)U+r}
=\sum_{r=1}^U {\mathcal B}_{jr} y_{(i-1)U+r} =x_{(i-1)l+j},
$$
for  each $i\in[n]$ and $j\in[l]$.
Hence, $\diag(X)=x$.
Moreover,
$$
A(I_n\otimes\mathfrak{w}^T)x = A (I_n\otimes\mathfrak{w}^T)Dy
= A(I_n\otimes\mathfrak{u}^T)y = b.
$$
Since $Y\succeq0$ and $\diag(Y)=y$, we have $y\geq0$.
Therefore, for each $i\in[n]$ we have
$$
\sum_{r=1}^U r y_{(i-1)U+r} \leq U\sum_{r=1}^U y_{(i-1)U+r} \leq U.
$$
Thus, $(I_n\otimes\mathfrak{w}^T)x = (I_n\otimes\mathfrak{u}^T)y  \leq Ue$.

Finally, using $D^T(I_n\otimes\mathfrak{w})
=I_n\otimes\mathfrak{u}$, we obtain
$$
\begin{aligned}
\left\langle (I_n \otimes \mathfrak{w})Q(I_n\otimes \mathfrak{w}^T),X \right \rangle
= \left\langle  D^T(I_n \otimes \mathfrak{w})Q (I_n\otimes \mathfrak{w}^T)D,Y \right\rangle
= \left\langle ( I_n \otimes \mathfrak{u})Q(I_n\otimes \mathfrak{u}^T), Y \right\rangle.
\end{aligned}
$$
\end{proof}

\medskip

Next, we compare the SDP relaxation \eqref{BSDPrelax} and the SDP relaxation obtained from the TSDP \eqref{ternarySDP_formulation_of_pure_integer_QP}.
After relaxing integrality constraints, the TSDP \eqref{ternarySDP_formulation_of_pure_integer_QP}
results in  the following SDP relaxation
\begin{equation}
\begin{array}{cl}  \label{TSDPrelax}
\min &
\left\langle
(I_n\otimes\mathfrak{v})Q (I_n\otimes\mathfrak{v}^{T}),X \right\rangle
+ (3^k-1)e^{T}Q (I_n\otimes\mathfrak{v}^{T})x + \frac{(3^k-1)^2}{4}e^{T}Qe \\[1ex]
\mathrm{s.t.}  & \displaystyle A(I_n\otimes\mathfrak{v}^{T})x = b-\frac{3^k-1}{2}Ae,\\[1ex]
& (I_n\otimes\mathfrak{v}^{T})x \leq \left( U-\frac{3^k-1}{2} \right)e \\[1ex]
&  \diag(X)\geq x, ~~ \diag(X)\geq -x,~~ \diag(X) \leq e \\[1ex]
& \begin{bmatrix}
1 & x^{T}\\
x & X
\end{bmatrix} \succeq 0,  \\[3ex]
\end{array}
\end{equation}
where $X \in {\mathcal S}^{nk}$  and $\mathfrak{v}$ defined in \eqref{mathfrakv}.
In the sequel we use the following partition of $X \in {\mathcal S}^{kn}$ into $n^2$ blocks, denoted by $X^{ij}$ for $i,j\in [n]$ of size $k\times k$, and similarly the decomposition $x^T = \begin{bmatrix} (x^1)^T, \cdots, (x^{n})^T \end{bmatrix}$ where $x^i \in \mathbb{R}^k$ for $i\in [n]$.

\begin{theorem}
Let $U \in \Z_+$ and $k= \left \lceil ~\log_3(U+1)~ \right \rceil$, $k\geq 1.$
 Then, the one-hot SDP relaxation~\eqref{BSDPrelax} dominates the ternary SDP relaxation~\eqref{TSDPrelax}.
\end{theorem}
\begin{proof}
Every $u\in[U]$ admits a unique   representation of the form~\eqref{sumx3}.
Let ${\mathcal T}\in\{-1,0,1\}^{k\times U}$ be the matrix such that
$\mathfrak{v}^T (\mathcal{T}+e_ke_U^{T}) = \mathfrak{u}^T.$
Define
$$
D:=\mathcal{T}+e_ke_U^{T}, \qquad
S:= I_n\otimes D \in\mathbb{R}^{nk\times nU}.
$$
Let $(y,Y)$  be feasible for~\eqref{BSDPrelax}.
We construct a feasible $(x,X)$ for~\eqref{TSDPrelax} as follows:
$$
\begin{bmatrix}
1 & x^{T}\\
x & X
\end{bmatrix}
:= \begin{bmatrix}
1 & 0\\
-e_{nk} & S
\end{bmatrix}
 \begin{bmatrix}
1 & y^{T}\\
y & Y
\end{bmatrix}
\begin{bmatrix}
1 & 0\\
-e_{nk} & S
\end{bmatrix}^T \succeq 0.
$$
A direct block multiplication gives
\begin{align}\label{defxx}
x =-e_{nk}+Sy
\end{align}
and
\begin{align}\label{defX}
X = J_{nk} -e_{nk}y^TS^T -Sye_{nk}^T +SYS^T.
\end{align}

Let us verify that ${X}$ satisfies the diagonal constraints from the SDP relaxation~\eqref{TSDPrelax}.
It follows from \eqref{defxx} that $x^i=-e_k+Dy^i$  for all $i\in[n]$.
For every $i\in [n]$, define
$$
q_i := 1-e_U^{T}y^i \geq 0,
$$
where the nonnegativity of $q_i$ follows from  $(I_n \otimes e^T_U) y \leq e $.
Using $D=\mathcal{T}+e_ke_U^{T}$  we obtain
\begin{align*}
x^i = -e_k +\left(\mathcal{T}+e_ke_U^T\right)y^i
= -e_k + \mathcal{T}y^i + e_k(e_U^Ty^i)
= -q_ie_k+\mathcal{T}y^i.
\end{align*}
It follows from \eqref{defX} that
$X^{ii}  = e_k e_k^T  -e_k (y^i)^T D^T -D y^i e_k^T +DY^{ii}D^T.$
From this and $Y^{ii}=\Diag(y^i)$ we obtain:
\begin{align*}
\diag(X^{ii}) & = e_k-2 Dy^i +\diag\left(D\Diag(y^i)D^T\right)
= e_k-2Dy^i+(D\circ D)y^i\\
& =  e_k+ \left(\mathcal{T}\circ\mathcal{T} -e_ke_U^T\right)y^i
=  \left ( 1-e_U^Ty^i \right )  e_k + (\mathcal{T}\circ\mathcal{T})y^i
=q_ie_k + (\mathcal{T} \circ \mathcal{T})y^i,
\end{align*}
where we used that $D\circ D-2D = \mathcal{T}\circ\mathcal{T}-e_ke_U^T.$
Therefore,
$$
\diag(X^{ii})- x^i
 = 2q_ie_k+ \left( \mathcal{T}\circ\mathcal{T}-\mathcal{T} \right)y^i \geq 0, \quad
\diag(X^{ii}) + x^i  = \left( \mathcal{T}\circ\mathcal{T}+\mathcal{T} \right)y^i \geq 0,
$$
and
$e_k-\diag(X^{ii}) = \left(  e_ke_U^T-\mathcal{T}\circ\mathcal{T} \right) y^i \geq 0.$
Thus,
$$
\diag(X^{ii}) \geq  x^i, \quad \diag(X^{ii}) \geq  -x^i, \quad  \diag(X^{ii}) \leq e_k, \quad \forall i\in [n].
$$
It remains to verify the equality and inequality constraints in~\eqref{TSDPrelax}.
By exploiting \eqref{defxx}, we obtain
$$
(I_n\otimes\mathfrak{v}^{T})x
=-(I_n\otimes\mathfrak{v}^{T})e_{nk}+
(I_n\otimes\mathfrak{u}^{T})y = -\frac{3^k-1}{2} e_n + (I_n\otimes\mathfrak{u}^{T})y,
$$
and thus
$$
A(I_n\otimes\mathfrak{v}^{T})x
=-\frac{3^k-1}{2}Ae_n+A (I_n\otimes\mathfrak{u}^{T})y
=-\frac{3^k-1}{2}Ae_n +b,
$$
where we used  $A(I_n\otimes\mathfrak{u}^{T})y=b$.
Hence, the equality constraint in~\eqref{TSDPrelax} is satisfied.

It follows from $(I_n\otimes e_U^{T})y\leq e_n$ that  $e_U^{T}y^i\leq1$ for $i\in[n]$.
Since $\mathfrak{u}\leq Ue_U$ and $y^i\geq0$, we obtain
$$
\mathfrak{u}^{T}y^i
\leq Ue_U^{T}y^i
\leq U,
\qquad i\in[n].
$$
Therefore, $(I_n \otimes \mathfrak{u}^{T})y \leq U e_n.$
Consequently,
\begin{align*}
(I_n\otimes\mathfrak{v}^{T})x =
(I_n\otimes\mathfrak{u}^{T})y -\frac{3^k-1}{2}e_n \leq
\left( U-\frac{3^k-1}{2} \right)e_n.
\end{align*}

It remains to verify that the objective values coincide.
Define  $V:= I_n\otimes\mathfrak{v}^{T}$ and
$L:= \begin{bmatrix} \frac{3^k-1}{2} e_n, ~ V \end{bmatrix}$.
Then
\begin{align*}
L
\begin{bmatrix}
1 & x^{T}\\
x & X
\end{bmatrix}
L^{T}
& =
L\begin{bmatrix}
1 & 0\\
-e_{nk} & S
\end{bmatrix}
\begin{bmatrix}
1 & y^{T}\\
y & Y
\end{bmatrix}
\begin{bmatrix}
1 & 0\\
-e_{nk} & S
\end{bmatrix}^{T}L^{T}
= (I_n\otimes\mathfrak{u}^{T} )Y(I_n\otimes\mathfrak{u}^{T} )^{T},
\end{align*}
and thus two objectives coincide:
\begin{align*}
\left\langle (I_n\otimes\mathfrak{u}^{T} )^{T}Q (I_n\otimes\mathfrak{u}^{T} ),Y\right\rangle
&
=
\left\langle Q,(I_n\otimes\mathfrak{u}^{T} )Y(I_n\otimes\mathfrak{u}^{T} )^{T}\right\rangle
=
 \left\langle Q, L
\begin{bmatrix}
1 & x^{T}\\
x & X
\end{bmatrix}
L^{T} \right\rangle \\
& =
\left\langle
V^{T}QV,X
\right\rangle
+ (3^k-1)\,e_n^{T}QV x + \frac{(3^k-1)^2}{4} e_n^{T} Q e_n.
\end{align*}
\end{proof}
The following examples show that the SDP relaxation~\eqref{BSDPrelax} may provide strictly better bounds than the other two relaxations, and that the latter two are incomparable.
\begin{example} \label{ex: 5.3}
Let $n=2$, $U=3$, and $k=l=2$.
Consider the following IQP:
$$ \min\left\{  x_1^2+3x_2^2:
x_1+x_2=3,\quad
x_i\in\{0,1,2,3\},\quad i\in[2]\right\}.
$$
Its optimal value is $7$, attained at $(x_1,x_2)=(2,1)$.

The optimal value of the binary SDP relaxation \eqref{BinExpSDPrelax} is ${27}/{4}$
the optimal value of the ternary SDP relaxation \eqref{TSDPrelax} is ${83}/{12}$, while the
optimal value of \eqref{BSDPrelax} is 7.
\end{example}

\begin{example} \label{ex: 5.4}
Let $n=2$, $U=3$, and $k=l=2$.
Consider the following IQP:
$$ \min\left\{ 2x_1^2+3x_2^2: x_1+x_2=4,\quad x_i\in\{0,1,2,3\},\quad i\in[2] \right\}.$$
Its optimal value is $20$, attained at $(x_1,x_2)=(2,2)$.

The optimal value of the ternary SDP relaxation \eqref{TSDPrelax} is ${59}/{3}$,
the optimal value of the binary SDP relaxation \eqref{BinExpSDPrelax} ${159}/{8}$,
while the optimal value of \eqref{BSDPrelax} is 20.
\end{example}
All of the SDP relaxations derived here may be  strengthened by adding valid inequalities, see e.g.,~\cite{Burer2014UnboundedCS,Galli2011,GALLI2021100661,McCormick1976}. However, we leave this for future work.

\section{MISDP for MIQP}
\label{sect: mixed-integer}
Next, we consider mixed-integer quadratic problems of the following form:

\begin{equation} \label{mixed_integer_QP}
   \begin{array}{cl}
      \min & \overline{x}^T Q \overline{x} \\
      \st    & A\overline{x} = b \\
              & \overline{x} \geq 0 \\
              & \overline{x}_i \in \Z \cap [ 0, U ] \qquad i \in \mathcal{B},
   \end{array}
\end{equation}
where $Q \in {\mathcal S}^n(\R)$, $A \in \R^{ m \times n }$, $b \in \R^m$, $U \in \Z_+$, and $\mathcal{B} \subset \{ 1,\ldots,n\}$ is the index set of the integer variables. W.l.o.g.\ we may assume that $\mathcal{B} = \{ 1,\ldots,k \}$.
After renaming the variables we rewrite~\eqref{mixed_integer_QP} as follows:
\begin{equation} \label{mixed_integer_QP_reformulated}
   \begin{array}{cl}
      \min & \begin{bmatrix} x \\ z \end{bmatrix}^T
                  \begin{bmatrix} Q_{11} & Q_{12} \\[1ex]
                  Q_{12}^T & Q_{22}\end{bmatrix}
                  \begin{bmatrix} x \\ z \end{bmatrix}  \\[2ex]
      \st    & A_1 x + A_2 z = b \\
       & z \geq 0 \\
              & x \in \Z^k \cap [ 0, U ]^k,
   \end{array}
\end{equation}
with $x \in \Z^k$, $z \in \R^{n-k}$, $Q_{11} \in \R^{k \times k}$,  $Q_{12} \in \R^{k \times n-k}$, $Q_{22} \in \R^{n-k \times n-k}$, $A_1 \in \R^{ m \times k }$ and $A_2 \in \R^{ m \times n-k }$.
As before, we define $\mathfrak{u} := \begin{bmatrix} 1, 2, \ldots, U \end{bmatrix} \in \Z^U$ and introduce $kU$ binary variables as follows:
$$
v^i_j := \begin{cases}
                0 & \text{ if } x_i \neq j \\
                1 & \text{ if } x_i = j,
            \end{cases}
$$
where $i \in [k]$ and $j \in [U]$.
Again, this gives rise to a matrix $V \in \R^{k \times U}$ defined by $V_{ij} := v^i_j$ such that $x = V \mathfrak{u}$, or equivalently $x = (I_k \otimes \mathfrak{u}^T) \vect (V^T).$
To ensure that this representation of $x$ is unique, we again impose the constraint
$ Ve \leq e$, see~\eqref{eq:Ve leq e}.

With this, we can can rewrite \eqref{mixed_integer_QP_reformulated} as follows:
\begin{equation} \label{mixedInteger as mixedBinary1}
   \begin{array}{cl}
      \min & \vect (V^T)^T (I_k \otimes \mathfrak{u}^T)^T Q_{11} (I_k \otimes \mathfrak{u}^T) \vect (V^T) + 2 z^T Q_{12}^{T} (I_k \otimes \mathfrak{u}^T) \vect (V^T) \\[1ex]
      & + z^T Q_{22} z  \\[1ex]
      \st   & A_1 (I_k \otimes \mathfrak{u}^T) \vect (V^T) + A_2 z = b \\[1ex]
            & (I_k \otimes e^T) \vect (V^T) \leq e \\[1ex]
            & \vect (V^T) \in \{ 0,1 \}^{kU} \\[1ex]
            & z \geq 0.
   \end{array}
\end{equation}
We assume in the following that  $z \leq  e U_R$  and rewrite \eqref{mixedInteger as mixedBinary1} as a mixed-binary SDP.
We first introduce the binary vector $y:= \vect (V^T)$.
It is well known that the linearization of the product of a binary variable and a nonnegative real variable is exact.
In particular,  one can replace the product of variables $W_{ij}:=y_i z_j$ for $i\in [kU]$, $j\in [n-k]$ with the following linear constraints
$$
W_{ij} \leq y_i U_R, \quad W_{ij} \leq z_j, \quad W_{ij} \geq z_j - U_R(1-y_i), \quad W_{ij}\geq 0.
$$
By exploiting these constraints, we can rewrite \eqref{mixedInteger as mixedBinary1} as follows:
$$
   \begin{array}{cl}
      \min & \langle (I_k \otimes \mathfrak{u}^T)^T Q_{11} (I_k \otimes \mathfrak{u}^T), Y \rangle
      + 2 \langle (I_k \otimes \mathfrak{u}) Q_{12}, W \rangle + \langle Q_{22}, Z \rangle  \\[1ex]
      \st   & A_1 (I_k \otimes \mathfrak{u}^T) y + A_2 {z} = b \\[1ex]
            & (I_k \otimes e^T) y \leq e \\[1ex]
            & W_{ij} \leq y_i U_R, \, W_{ij} \leq z_j, \, W_{ij} \geq z_j - U_R(1-y_i) \quad  i\in [kU],\, j\in [n-k]\\[1ex]
            & \begin{bmatrix} 1 & y^T \\ y & Y \end{bmatrix} \succeq 0, ~~ \diag (Y) = y \\[2ex]
            &  \begin{bmatrix} 1 & z^T \\ z & Z \end{bmatrix} \succeq 0, ~~  {\rm rank}\begin{bmatrix} 1 & z^T \\ z & Z \end{bmatrix}=1 \\[2ex]
            & {z} \geq 0, \ W \geq 0, \, y \in \{ 0,1 \}^{kU}.
   \end{array}
$$
Note that the above optimization problem does not follow  into our definition of
an MISDP problem \eqref{MISDP}  due to the rank-one constraint.
Nevertheless, after relaxing integrality and the rank-one constraint one can obtain an SDP relaxation for the MIQP \eqref{mixed_integer_QP}.
The resulting SDP relaxation  is, to the best of our knowledge, the only SDP relaxation of a general mixed-integer QP in the literature.

\section*{Declaration of AI use}
Although the authors resisted using AI to write the manuscript until the final stage, they ultimately gave in.
We used ChatGPT to assist with reviewing the manuscript and refining its exposition, which sped up the final process. AI was also used to find Example \ref{ex: 5.3} and Example \ref{ex: 5.4}.
The authors reviewed all AI-assisted material and take full responsibility for the manuscript.

\bibliographystyle{plain}
\bibliography{DuerSotirov}

\end{document}